%% file: ms.tex
\documentclass[11pt]{article}
\usepackage{amsmath,amssymb,amsthm,bm}
\usepackage{mathtools}
\usepackage{dsfont}
\usepackage{verbatim}
\usepackage{graphicx,rotating}
\usepackage{xcolor}
\usepackage{ulem} % crossing out text
\usepackage{enumitem} % shift itemize environment
\usepackage{xr}
\usepackage[colorlinks,citecolor=blue,urlcolor=blue]{hyperref}

\input{\auxpath self_shortcuts.tex}

\input{\auxpath self_environments.tex}

\newcommand{\cl}{{\rm cl}}

\newcommand{\pdM}{{\partial_d \cM}}

\begin{document}

\title{ On the finiteness of the second moment of the number of critical points of Gaussian random fields }
\author{Samuel Davenport$^1$, Fabian J.E. Telschow$^{2, }$\footnote{Both authors contributed equally.} \\[5mm]
	$^1$Division of Biostatistics, University of California, San Diego \\
	$^2$Institute of Mathematics, Humboldt Universit\"at zu Berlin
}
\date{\today}
\maketitle

%========================================================================
% Abstract+
%========================================================================
\input{\auxpath abstract.tex}

%========================================================================
% Introduction
%========================================================================
\input{\auxpath introduction.tex}

%========================================================================
% Main document
%========================================================================
\input{\auxpath main_body.tex}

%========================================================================
% Acknowledgements
%========================================================================
\section*{Acknowledgments}
F.T. is funded by the Deutsche Forschungsgemeinschaft (DFG) under Excellence Strategy The Berlin Mathematics Research Center MATH+ (EXC-2046/1, project ID:390685689). S.D. was partially supported by NIH grant R01EB026859.

%========================================================================
% Bibliography
%========================================================================
\bibliographystyle{plain}

%\bibliography{\auxpath paper-ref.bib,\auxpath RFT.bib}

%========================================================================
% Appendix
%========================================================================
\input{\auxpath appendix.tex}

\end{document}

%% file: self_shortcuts.tex
% Mathbb

\newcommand{\mE}{\mathbb{E}}

% Mathcal

\newcommand{\cM}{\mathcal{M}}

% Matrices

% Stochastics

\newcommand{\Cov}{{\rm Cov}}

% RFT

%% file: self_environments.tex
\newtheorem{theorem}{Theorem}%[section]
\newtheorem{lemma}{Lemma}%[section]
%[section]
%[section]

\theoremstyle{definition}
\newtheorem{definition}{Definition}%[section]
%[section]
\theoremstyle{plain}
\newtheorem{assumption}{Assumption}

\theoremstyle{remark}

\newtheorem{remark}{Remark}

%% file: abstract.tex
\begin{abstract}
	We prove that the second moment of the number of critical points of any sufficiently regular random field, for example with almost surely $ C^3 $ sample paths, defined over a compact Whitney stratified manifold is finite. Our results hold without the assumption of stationarity - which has traditionally been assumed in other work. Under stationarity we demonstrate that our imposed conditions imply the generalized Geman condition of \cite{Estrade2016}.
\end{abstract}

%% file: introduction.tex
\section{Introduction}
The variance of the number of critical points of a Gaussian random field $ f $ has been studied by many authors over the last 60 years. Papers considering this topic have typically assumed that the random field is stationary and that the underlying domain of interest is a compact subset of $ \mathbb{R}^D $ where $ D \in \mathbb{N}$ denotes the dimension of the random field. In our work we generalize these results to non-stationary fields on arbitrary $ C^3 $ manifolds, showing that, under certain regularity conditions, the variance of the number of critical points of the random field is finite. 

Formulae for the moments of the expected number of critical points were originally developed by \cite{Rice1944}. These were generalized to formulae for the expected factorial $ n $th moment of the number of critical points (for $ n \in \mathbb{N} $) in \cite{Cramer1965}, \cite{Belyaev1966} and \cite{Belyaev1967b} in 1D, see also \cite{Cuzick1975}. \cite{Malevich1989} extended these results to fields on domains of arbitrary dimension, see \cite{Adler2010} for a comprehensive overview of the proof. These papers showed that it is sufficient to show that the Rice integral for the $ n $th factorial moment converges in order to show that the $ n $th moment of the number of critical points is finite.

Using this approach, applied to the second moment and assuming stationarity, simplifies the problem and has lead to a number of results. Working in 1D \cite{Cramer1967} showed, letting $ r $ denote the covariance function of the random field, that if there exists some $ \delta > 0 $ such that 
\begin{equation*}
\int_{0}^{\delta} \frac{r^{(4)}(0) - r^{(4)}(t)}{t} \,dt
\end{equation*}
is finite then this is sufficient (where $ r^{(4)} $ denotes the fourth derivative of $ r $). \cite{Geman1972} showed that this is also a necessary condition. It was only recently that \cite{Estrade2016} extended these results to stationary Gaussian random fields on $\mathbb{R}^D, D \in \mathbb{N}$, showing that (letting $ r $ denote the now multivariate covariance function) if there exists some $ \delta > 0 $ such that 
\begin{equation*}
\int_{\left\lVert t \right\rVert < \delta} \frac{\left\lVert r^{(4)}(0) - r^{(4)}(t) \right\rVert}{\left\lVert t \right\rVert^D} \,dt
\end{equation*}
is finite (where $r^{(4)}(t)$ denotes the 4th derivative), then this generalized Geman condition is sufficient. A similar result was stated in \cite{Elizarov1985}, however the proof presented there lacks sufficient detail to make it rigorous. Under the additional assumption of isotropy, \cite{Estrade2016b} showed that the result holds under minimal other assumptions.

When the assumption of stationarity is removed the problem becomes more difficult. The reason for this is that the Rice integral integrates twice over the same domain and this leads to non-degeneracies within the integrand which are difficult to bound. Recently, by adapting the arguments of \cite{Piterbarg1996}, \cite{Cheng:2015extremes} showed that, under regularity assumptions, the second factorial moment of the number of critical points within an epsilon ball around a given point is $ o\big( \epsilon^D \big) $. In this paper we use this result to show that the second moment of the number of critical points of  Gaussian random fields over certain classes of compact stratified spaces is finite. The classes of stratified spaces we consider include compact manifolds with and without boundary and Whitney-stratified spaces (e.g., \cite[Chapter 8.1]{Adler:RFT2009}. 

%combined with compactness to show that each given manifold can be broken down into finitely many such balls and that, as a result, the second moment of the number of critical points within the manifold as a whole is finite. 

This paper is laid out as follows. Section 2 introduces the notation and definitions that we will need and then lays out the assumptions under which our results hold. In Section 3 we state and prove a lemma that shows that our assumptions hold under diffeomorphic transformations and include our main result regarding the finiteness of the second moment.

%% file: main_body.tex
\section{Notation and Definitions}
In this section we shall introduce a number of the definitions and key assumptions upon which our results depend. Let $ (\, \tilde{\mathcal{M}}, g \,) $ be a $ D $-dimensional $ C^3 $-Riemannian manifold without boundary for some $D \in \mathbb{N}$, where $\mathbb{N}$ denotes the set of positive integers. Throughout the article any submanifold $\mathcal{M}\subset \tilde{\mathcal{M}}$ is assumed to be an embedded submanifold, i.e. such that the map $\iota:~\mathcal{M} \rightarrow \tilde{\mathcal{M}}$, $x\mapsto x$, is an embedding. We shall denote the dimension of a (sub)-manifold $\mathcal{N}$ by $\text{dim}(\mathcal{N}) \in \{\, 1,\ldots, D \,\,\}$. 
Recall that a compact manifold without boundary is called \emph{closed}.
For any $ \mathcal{A} \subset \tilde{\mathcal{M}} $ we shall denote the toplogical closure/interior of $\mathcal{A}$ as $\text{cl}(A) / \text{int}(A)$ respectively.

In order for our results to hold we shall require the following notion of extendability.
\begin{definition}
	We say that a submanifold $ \mathcal{N} \subset \tilde{\mathcal{M}} $ is extendable, if there exists a
	$ \text{dim}(\mathcal{N}) $-dimensional, $ C^3 $-submanifold $ \mathcal{A} \subset \tilde{\mathcal{M}} $, $ \text{dim}(\mathcal{N}) \in \{\, 1, \ldots, D \,\} $, without boundary, such that $ \cl\big( \mathcal{N} \big) $ is a $ C^3 $-submanifold of $ \mathcal{A} $.
\end{definition}
This extendability condition is widely satisfied. First, $ \tilde{\mathcal{M}} $ itself is trivially extendable.
Secondly, if $ \cl\big( \mathcal{ N } \big) $ is a $ D $-dimensional compact
$ C^3 $-submanifold with boundary $ \partial \mathcal{ N } $ such that
$  \partial \mathcal{ N } $ is a $ (D - 1) $-dimensional orientable $ C^3 $-submanifold
of $ \tilde{\mathcal{M}} $, then $ \mathcal{N} $ is extendable, since $ \cl\big( \mathcal{ N } \big) $ can be embedded into the union of $ \cl\big( \mathcal{ N } \big) $ and a tubular neighbourhood, which exists by \cite[Exercise 8-5.]{Lee1997}.
Thirdly, any open $d$-dimensional polygon $\mathcal{P}$, $d\in \{\,1,\ldots,D\,\}$, where $ \cl\big( \mathcal{P} \big) $ is contained in an open subset $ \mathcal{O} $ of an $ d $-dimensional affine plane of $ \mathbb{R}^D $ is extendable, since $ \mathcal{O} $ can be used for $ \mathcal{A} $. In particular this implies that the interior, all faces and all edges of the unit cube
$ [-1,1]^D \subset \mathbb{R}^D $ are extendable. 
Moreover the same applies to many curved submanifolds which are themselves proper $ d $-dimensional submanifolds of a $ d $-dimensional submanifold of $\tilde{\mathcal{M}}$, for example, $ \mathbb{S}^{D-1} =  \{\, x \in \mathbb{R}^{D}:~ x_1^2 + \ldots + x_D^2 = 1\,,~ x_D > c \,\} $ for some $c \in (-1, 1)$.
Finally, if $ \cl\big( \mathcal{ N } \big) $ is a $ d $-dimensional closed $ C^3 $-submanifold, then $\mathcal{ N }$ is extendable. This implies that spaces like the unit sphere $\mathbb{S}^{D-1}$ and
$ \mathbb{S}^{D-1} \backslash \{\,(1,0,\ldots,0)\,\} $ are extendable.

To define critical values of a function $ h \in C^1( \mathcal{N} ) $ recall
that any vector field $ X:~\mathcal{N} \rightarrow T\mathcal{N} $, where
$ T\mathcal{N} $ denotes the tangent bundle of $ \mathcal{N} $, defines
a function $ X\!h: \mathcal{N} \rightarrow \mathbb{R}$.

\begin{definition}
	Let $\mathcal{N} $ be a $d$-dimensional submanifold of $ \tilde{\mathcal{M}} $, then the set of critical points of $ h\in C^1( \mathcal{N} ) $ is defined as
	\begin{equation}\label{eq:DefCrits}
	\mathcal{C}_h(\mathcal{N})	= \left\{ s \in \mathcal{N} :~
	X\!h (s) = 0 \text{ for all vector fields }
	X \text{ on } \mathcal{N}
	\,\right\}\,.
	\end{equation}
	We will denote the number of critical points by
	\begin{equation}\label{eq:NumCrits}
	\mu_h\big( \mathcal{N} \big)
	= \# \mathcal{C}_h\big( \mathcal{N} \big)\,.
	\end{equation}
\end{definition}
\begin{remark}
	Let $ \mathcal{U} \subset \mathcal{N} $ and $ \mathcal{V} \subset\mathbb{R}^{\text{dim}(\mathcal{N})} $ be open and
	$ \varphi: \mathcal{U} \rightarrow \mathcal{V} $ be a differentiable chart of $ \mathcal{N} $.
	For $ i \in \{\, 1,\ldots, {\text{dim}(\mathcal{N})} \,\} $ we define the vector field
	$ E_{i}: \mathcal{U} \rightarrow T \mathcal{U} $
	by
	$ E_i h(s) = \frac{d}{dt}\left(\, (h \circ \varphi^{-1})( \varphi(s) + t e_i ) \,\right)\big\vert_{t=0} $, where $ e_i $ is the $i$-th standard Euclidean basis vector.
	For all $s \in \mathcal{U}$, the tangent vectors $ E_1(s), \ldots, E_{\text{dim}(\mathcal{N})}(s) $ form a basis of the tangent space
	$ T_s \mathcal{U} $ and so it follows that $ E_ih (s) = 0 $ for all $ i\in\{\, 1, \ldots, {\text{dim}(\mathcal{N})} \,\} $
	if and only if $ X\!h(s) = 0 $ for all vector fields $ X $.
	Moreover, $ E_ih (s) = 0 $ if and only if
	$ \frac{\partial (h\circ\varphi^{-1}) }{ \partial x_i}\big(\varphi(s)\big) = 0 $.
	Therefore $ s \in \mathcal{N} $ is a critical point in the manifold sense defined by  $ X\!h(s) = 0 $
	for all vector fields $ X $	if and only if there exists a chart $ \varphi $ such that $ \varphi(s) \in \mathbb{R}^{\text{dim}(\mathcal{N})} $
	is a critical point of  $ h\circ \varphi^{-1} $ in the usual Euclidean sense.
\end{remark}

The set of critical values and number of critical values can be naturally extended to functions $ h\in C^1( \mathcal{M} )$ defined on disjoint unions
$ \mathcal{M} = \bigsqcup_{ d = 0 }^D \pdM \subset \tilde{\mathcal{M}} $, where each $ \pdM $
is either empty or a $d$-dimensional $ C^3 $-submanifold without boundary.
Letting $ h\vert_\pdM $ denote the restriction of $h$ to $ \pdM $, this is done by setting
\begin{equation}\label{eq:DefCritsII}
\mathcal{C}_h(\mathcal{M})	= \bigcup_{ d = 0 }^D
\mathcal{C}_{h\vert_\pdM}(\pdM)\,,~ ~ ~  
\mu_h\big( \mathcal{M} \big) = \sum_{d=0}^D		\mu_{h\vert_\pdM}\big( \pdM \big)\,.
\end{equation}
This broader definition means that our results will apply, among others, to manifolds with boundary, manifolds with corners and more general Whitney stratified manifolds (e.g., \cite[Chapter 8.1]{Adler:RFT2009}).

In particular note that, if $ \mathcal{M} $ is a Whitney stratified manifold, the set of critical points defined in \eqref{eq:DefCritsII} is the same as the set of critical points defined in \cite[p. 194]{Adler:RFT2009}. To see this recall that if $\tilde{\mathcal{M}} $ has a Riemannian metric $ g $,
then $ g $ induces a Riemannian metric on any $ C^3 $-submanifold $ \mathcal{N} \subset \tilde{\mathcal{M}} $ without boundary.
We will also denote this induced Riemannian metric by $ g $.
The gradient of $ h \in C^1(\, \mathcal{N} \,) $, denoted by  $ \nabla^{\mathcal{N}}\! h $, is the unique continuous vector field on $ \mathcal{N} $ such that $ g\big(\, \nabla^{\mathcal{N}}\!h, X \,\big) = X\!h $ for every vector field $ X $ on $ \mathcal{N} $. Hence given $s \in \mathcal{N}$, $ \nabla^{\mathcal{N}}\! h(s) = 0 $ if and only if $ X\!h(s) = 0 $ for every vector field $ X $ on $ \mathcal{N} $, which establishes the
connection between the two definitions.

Recall that $ (\,\mathcal{U}, \varphi \,) $ is a $C^3$-chart around $ s \in \tilde{\mathcal{M}} $ if
$\mathcal{U}\subset \tilde{\mathcal{M}}$ is open and contains $ s $ and $ \varphi\in C^3(\,\mathcal{U}, \mathcal{V}\,) $ is a diffeomorphism onto an open set $ \mathcal{V} \subset \mathbb{R}^D $. For such a chart and a function
$ h \in C^2( \tilde{\mathcal{M}} )$ we define the functions $ h^{\varphi} = h \circ \varphi^{-1} $ and its derivatives as
$ h^{\varphi}_i =  \partial h^{\varphi} / \partial x_i  $ and
$ h^{\varphi}_{ij} =  \partial^2 h^{\varphi} / \partial x_i\partial x_j  $ for $i,j=1,\ldots,D$.
Moreover, we write
$ \nabla h^{\varphi} = \big(\, h^{\varphi}_1, \ldots, h^{\varphi}_D \,\big) $ for the gradient of $ h^{\varphi} $ and
$ \nabla^2{h^{\varphi}} $ for the Hessian of $ h^{\varphi} $, which has as
the $(i,j)$-th entry $ h^{\varphi}_{ij} $.

Furthermore, given $D, D' \in \mathbb{N}$ and a matrix $A \in \mathbb{R}^{D \times D'}$ we define its \emph{vectorization} as $ \mathbf{vec}( A ) = (\, A_{11},\ldots, A_{D'1}, \ldots, A_{1D},\ldots, A_{DD'}  \,)^T $, and the \emph{half-vectorization}  of a symmetric matrix $A \in \mathbb{R}^{D \times D}$ as $ \mathbf{vech}( A ) = (\, A_{11}, \ldots, A_{D1}, A_{22}, \ldots, A_{D2}, \ldots, A_{D-1D-1}, A_{DD-1}, A_{DD}\,)^T $.

With these definitions we formulate the following assumption on a random field $ f $ on $ \tilde{\mathcal{M}} $.
\begin{assumption}\label{assumption:Fields}
	$ f $ has almost surely $ C^2( \tilde{\mathcal{M}} ) $
	sample paths and for all $s \in \tilde{\mathcal{M}}$ there exists
	a $C^3$-chart $ (\,\mathcal{U}, \varphi \,) $ around
	$ s $ mapping $\mathcal{U} \subset \tilde{\mathcal{M}}$ to an open set $\mathcal{V}\subset\mathbb{R}^D$, such that the following conditions hold.
	\begin{itemize} 
		\item[\textbf{(C1)}] 
		There exist positive constants $ L \in (0, \infty ) $ and $\eta > 0 $ such that
		$$
		\mathbb{E}\big[\, f_{ij}^{\varphi}(x) - f_{ij}^{\varphi}(y) \,\big]^2 \leq L \Vert\, x - y \,\Vert^{2\eta}\,,
		$$
		for all $ x, y \in  \mathcal{V} $ and all $ i,j \in \{\, 1, \ldots, D \,\} $. 
		\item[\textbf{(C2)}]
		For each $ (x, y) \in \mathcal{V} \times \mathcal{V} $
		with $x \neq y$, the Gaussian random vector
		$$ \Big(\, %f^{\varphi}(x),
		\nabla f^{\varphi}(x),
		\mathbf{vech}\big( \nabla^2{f^{\varphi}}(x) \big)^T,
		%f^{\varphi}(y),
		\nabla f^{\varphi}(y),
		\mathbf{vech}\big( \nabla^2{f^{\varphi}}(y) \big)^T \,\Big)^T $$
		is non-degenerate in the sense that its covariance matrix is positive definite.
	\end{itemize}
\end{assumption}
Given a $C^3$-chart $(\,\mathcal{U}, \varphi\,)$ such that \textbf{(C1)}/\textbf{(C2)} holds with respect to a random field $f$ on $\tilde{\mathcal{M}}$, we shall say that $f^{\varphi}$ satisfies \textbf{(C1)}/\textbf{(C2)} on $\mathcal{U}$ and when the existence of $f$ is clear we will simply say that $(\,\mathcal{U}, \varphi\,)$ satisfies \textbf{(C1)}/\textbf{(C2)} and take $f$ to be implicit.

\begin{remark}
	\cite{Cheng:2015extremes} introduced slightly different conditions (their (C1') and (C2')) in order to obtain a manifold version of their Lemma 4.2.
	These are not quite what we want, since
	if they are assumed for a random field $ f $
	over $ \tilde{\mathcal{M}} $, they do not necessarily hold for the restriction of $ f $ to arbitrary submanifolds $\mathcal{N}$ of $\tilde{\mathcal{M}}$. This is because for an arbitrary submanifold $ \mathcal{N} \subset \tilde{\mathcal{M}} $ in general there is no collection of $\text{dim}(\mathcal{N})$ orthonormal fields of the $D$ orthonormal fields chosen on $ \tilde{\mathcal{M}} $ such that all of these fields are orthonormal fields for $\mathcal{N}$. However, in Lemma \ref{lem:TransformationAssumption} we show that, if \textbf{(C1)}, \textbf{(C2)} hold in any chart containing $s \in \tilde{\mathcal{M}}$, then they hold in all charts containing $s$. Thus, in particular they hold for the the special charts chosen in \cite{Cheng:2015extremes}. Note also that the inclusion of $f(x)$ and $f(y)$ in the vector in \textbf{(C2)} is unnecessary, despite its inclusion in the conditions of \cite{Cheng:2015extremes}, as it was not used in the proof of their Lemma 4.2.
\end{remark}

\begin{remark}
	If $\tilde{\mathcal{M}} \subset \mathbb{R}^D$ and $f$ is stationary and satisfies \textbf{(C1)} then it satisfies the generalised Geman condition, see Lemma \ref{lem:Geman}.
\end{remark}

%========================================================================
% Proof of main result
%========================================================================
\section{Proof of the main result}

In order to prove our main result we first prove an important technical lemma. This lemma justifies the use of Conditions \textbf{(C1)} and \textbf{(C2)} since it shows that they are compatible with changing coordinates. We will use this, in Lemma \ref{lem:SimpleSubmanifolds}, to show that if these conditions hold for a random field $ f $ on $ \tilde{\mathcal{M}} $, then they also hold for the restriction of $ f $ to any $ C^3 $-submanifold $ \mathcal{N} $ of $ \tilde{\mathcal{M}} $.

\begin{lemma}\label{lem:TransformationAssumption}
	Let $ f $ be a Gaussian random field on
	$ \tilde{ \mathcal{M} } $ satisfying Assumption
	\ref{assumption:Fields} for the $C^3$ chart $ (\,\mathcal{U}, \varphi\,) $.
	Let $ (\,\mathcal{U}', \phi\,) $ be another
	$C^3$-chart such that $ \mathcal{W} = \mathcal{U} \cap \mathcal{U}' $ is non-empty.
	Then for any $s \in \mathcal{W}$ there exists an open set $ s \in \mathcal{W}_s \subset \mathcal{W} $ such that $ (\,\mathcal{W}_s, \phi\,) $ satisfies Assumption \ref{assumption:Fields}.
\end{lemma}  
\begin{proof}
	Since $ f^{\phi} = f^{\varphi} \circ \varphi \circ
	\phi^{-1} $ on $ \phi(\mathcal{W}) $, the chain rule implies that for $ \psi = (\, \psi_1,\ldots, \psi_D\,) = \varphi \circ \phi^{-1},$
	\begin{equation}\label{eq:1stDeriv}
	f^\phi_{i} = \big( f^\varphi\circ\psi\big)_{i}
	= \sum_{d=1}^D \frac{\partial f^\varphi\circ\psi}{ \partial \psi_d } \frac{\partial \psi_d}{\partial s_i}
	= \sum_{d=1}^D f^\varphi_d\circ\psi \frac{\partial \psi_d}{\partial s_i}\,,~ ~ i = 1,\ldots, D\,,
	\end{equation}
	which shows that $ \nabla f^\phi(x) = \nabla f^{\varphi}(\tilde x ) J_\psi(x)$ for all $x\in \phi(\mathcal{W})$, where $ \tilde x = \psi(x)$ and $ J_\psi(x) $ is the Jacobian matrix of $ \psi $ at $x$. For the second order derivatives we obtain:
	\begin{equation}\label{eq:2ndDeriv}
	\begin{split}
	f^\phi_{ij}
	&= \sum_{d=1}^D \frac{\partial f^\varphi\circ\psi}{ \partial \psi_d } \frac{\partial^2 \psi_d}{\partial s_i \partial s_j} + \sum_{d,d'=1}^D \frac{\partial ^2f^\varphi\circ\psi}{ \partial \psi_d\partial \psi_{d'} } \frac{\partial \psi_d}{\partial s_i }\frac{\partial \psi_{d'}}{ \partial s_j }\\
	&= \sum_{d=1}^D f^\varphi_d\circ\psi \frac{\partial^2 \psi_d}{\partial s_i \partial s_j} + \sum_{d,d'=1}^D f^\varphi_{dd'}\circ\psi \frac{\partial \psi_d}{\partial s_i }\frac{\partial \psi_{d'}}{ \partial s_j }	\,,~ ~ i,j = 1,\ldots, D\,.
	\end{split}
	\end{equation}
	This can be written in vector notation as
	\begin{equation}
	\mathbf{vech}\Big(\, \nabla^2{f^{\phi}}(x) \,\Big)
	= L\,\mathbf{vec}\Bigg(\, \sum_{d=1}^D f_d^{\varphi}\big( \tilde x \big)\nabla^2{\psi_d}(x) \,\Bigg) + L\, \big(\, J_\psi(x) \otimes J_\psi(x) \,\big)\,R\,\mathbf{vech}\Bigg( \nabla^2{f^{\varphi}} \big( \tilde x \big) \Bigg)\,.
	\end{equation}
	Here $ L \in \mathbb{R}^{D(D+1)/2 \times D^2}
	$ is the elimination matrix
	and $ R \in \mathbb{R}^{D^2 \times D(D+1)/2} $ is the duplication matrix, the precise definitions of which can be found in 
	\cite{Magnus1980}. The matrix $L \big(\, J_\psi(x) \otimes J_\psi(x) \,\big)R$ is invertible by Lemma 4.4.iv of \cite{Magnus1980}, and the fact that $J_\psi(x)$ is invertible because $ \psi $ is a diffeomorphism.
	Furthermore note that, letting ${\rm detcov}( U ) = \text{det}\big(\,\Cov[U]\,\big)$, for any random vector $U \in \mathbb{R}^D$, 
	\begin{equation}\label{eq:detcov}
	{\rm detcov}( U ) = {\rm detcov}( MU ) \,/ \det(M)^{2D}
	\end{equation}
	for any invertible matrix $ M \in \mathbb{R}^{D \times D} $, see Lemma \ref{lem:MU}.
	
	\textit{Proof of} \textit{\textbf{(C2)} for the chart} $(\, \mathcal{W}, \phi \,)$\textit{:}  We shall show that for all $ x, y \in \phi(\mathcal{W})$ the distribution of the random vector $V$ defined by
	\begin{equation*}
	\begin{split}
	V &= \Big(\, 
	\nabla^T f^{\phi}(x);
	\mathbf{vech}\Big(\, \nabla^2{f^{\phi}}(x) \,\Big);
	\nabla^T f^{\phi}(y);
	\mathbf{vech}\Big(\, \nabla^2{f^{\phi}}(y) \,\Big) \,\Big)\\
	&=
	\Bigg( 
	J_\psi(x)^T\nabla^T f^{\varphi}(\tilde x);
	L\,\mathbf{vec}\Bigg(\, \sum_{d=1}^D f_d^{\varphi}\big( \tilde x \big)\nabla^2{\psi_d} \,\Bigg) + L\, \big(\, J_\psi(x) \otimes J_\psi(x) \,\big)\,R\,\mathbf{vech}\Big(\, \nabla^2{f^{\varphi}} \big( \tilde x \big) \,\Big);\\
	&~ ~ ~ ~ ~ ~ ~ ~ ~  ~ 
	J_\psi(y)^T\nabla^T f^{\varphi}(\tilde y);
	L\,\mathbf{vec}\Bigg(\, \sum_{d=1}^D f_d^{\varphi}\big( \tilde y \big)\nabla^2{\psi_d} \,\Bigg) + L\, \big(\, J_\psi(y) \otimes J_\psi(y) \,\big)\,R\,\mathbf{vech}\Big(\, \nabla^2{f^{\varphi}} \big( \tilde y \,\big) \,\Big) \Bigg),
	\end{split}
	\end{equation*}
	where $\tilde{y} = \psi(y)$ (and we have used ; to vertically stack vectors), is non-degenerate.
	Let $I_{D\times D}$ denote the identity matrix in $\mathbb{R}^{D\times D}$. Multiplying $V$ by the the block diagonal matrix $B$ which has diagonal blocks
	$$\big(J_\psi^T(x)\big)^{-1}, I_{D(D+1)/2 \times D(D+1)/2}, \big(J_\psi^T(y)\big)^{-1}, I_{D(D+1)/2 \times D(D+1)/2} $$
	only changes $\text{detcov}[V]$ by a positive scalar by applying  \eqref{eq:detcov}.
	Scaling the gradient appropriately and adding it to the vectors in
	the third and sixth components it follows that
	\begin{equation*}
	\begin{split}
	&\Bigg( 
	\nabla^T f^{\varphi}(\tilde x);
	L\, \big(\, J_\psi(x) \otimes J_\psi(x) \,\big)\,R\,\mathbf{vech}\Big(\, \nabla^2{f^{\varphi}} \big( \tilde x \big) \,\Big);\\
	&~~~~~~~~~~~~~~~~~~~~~
	\nabla^T f^{\varphi}(\tilde y);
	L \big(\, J_\psi(y) \otimes J_\psi(y) \,\big)\,R\,\mathbf{vech}\Big(\, \nabla^2{f^{\varphi}} \big( \tilde y \big) \Big) \,\Bigg)
	\end{split}
	\end{equation*}
	has, up to a positive scalar, the same $ {\rm detcov} $ as the random vector $V$.
	For all $x\in \phi(\mathcal{W}) $, $ L \big(\, J_\psi(x) \otimes J_\psi(x) \,\big)R $ is invertible, and so applying \eqref{eq:detcov} once more with an appropriate block diagonal matrix shows that, up to a positive scalar, the random vector
	\begin{equation*}
	\begin{split}
	&\Big(\, %f^{\varphi}(\tilde x),
	\nabla^T f^{\varphi}(\tilde x);
	\mathbf{vech}\big(\,\nabla^2{f^{\varphi}}(\tilde x)\,\big);
	%f^{\varphi}(\tilde y),
	\nabla^T f^{\varphi}(\tilde y);
	\mathbf{vech}\big(\,\nabla^2{f^{\varphi}}(\tilde y)\,\big) \,\Big)
	\end{split}
	\end{equation*}
	has the same $ {\rm detcov} $ as $V$. Since \textit{\textbf{(C2)}} holds for the chart $\varphi$ this last vector is non-degenerate for all $\tilde x, \tilde y \in \varphi(\mathcal{W})$ and hence $V$ is non-degenerate. 
	
	\textit{Proof of} \textit{\textbf{(C1)} for a chart} $(\, \mathcal{W}_s, \phi \,)$ \textit{for} $s\in\mathcal{W}$\textit{:} Since $ \mathcal{W} $ is non-empty and open there exists an open set $ \mathcal{W}_s $ around $ s $ such that $ \cl\left( \mathcal{W}_s \right) \subset \mathcal{W} $.
	
	Using equation \eqref{eq:2ndDeriv} and the fact that $ \big( \sum_{n=1}^N a_n \big)^2 \leq 2^{N-1} \sum_{n=1}^N a_n^2 $ for all $a_1,\ldots, a_N \in \mathbb{R}$, proving \textit{\textbf{(C1)}} reduces to showing that for all $d,d',i,j \in \{\, 1,\ldots, D \,\}$
	\begin{align*}
	\mathbb{E}\Bigg[\, \Bigg( f^\varphi_d\big( \psi(x) \big) \frac{\partial^2 \psi_d}{\partial s_i \partial s_j}(x) - f^\varphi_d\big( \psi(y) \big) \frac{\partial^2 \psi_d}{\partial s_i \partial s_j}(y) \Bigg)^2 \,\Bigg]
	\leq M \Vert\, x - y \,\Vert^{2\eta'}
	\end{align*}
	
	\begin{align*}
	\text{ and }\mathbb{E}\Bigg[\, \Bigg( f^\varphi_{dd'}\big( \psi(x) \big) \frac{\partial \psi_d}{\partial s_i }(x)\frac{\partial \psi_d}{ \partial s_j }(x) - 
	f^\varphi_{dd'}\big( \psi(y) \big) \frac{\partial \psi_d}{\partial s_i }(y)\frac{\partial \psi_{d'}}{ \partial s_j }(y) \Bigg)^2 \,\Bigg]
	\leq M \Vert\, x - y \,\Vert^{2\eta'}
	\end{align*}
	for some $ \eta' > 0 $, some $ M \in (0,\infty)$ and all $x,y\in \phi(\mathcal{W}_s)$. In order to bound the left hand side of these expressions, notice that given a function $g: \tilde{\mathcal{M}} \rightarrow \mathbb{R}$ and a random field $Y: \tilde{\mathcal{M}} \rightarrow \mathbb{R}$,
	\begin{align*}
	&\mathbb{E}\Bigg[\, \Bigg( g(x)Y\big(\psi(x)\big) - g(y)Y\big(\psi(y)\big) \Bigg)^2 \,\Bigg]\\
	&~~~~~~~~~~\leq 2g(x)^2\, \mathbb{E}\Bigg[ \,\Bigg( Y\big(\psi(x)\big) - Y\big(\psi(y)\big) \Bigg)^2 \,\Bigg]
	+2\mathbb{E}\big[\,Y^2\big(\psi(y)\big)\,\big]\big(\, g(x) - g(y) \,\big)^2\\
	&~~~~~~~~~~\leq 2 \max_{x \in \phi(\mathcal{W}_s)}\{\, g^2(x)\,\} \mathbb{E}\Bigg[\, \Bigg( Y\big(\psi(x)\big) - Y\big(\psi(y)\big) \Bigg)^2 \,\Bigg]
	+2\max_{x\in \varphi(\mathcal{W}_s)}\left\{\,\mathbb{E}\big[\,Y^2(x) \,\big] \,\right\} \big(\, g(x) - g(y) \,\big)^2.
	\end{align*}
	Taking $ g = \frac{\partial \psi_d}{ \partial s_i }\frac{\partial \psi_{d'}}{ \partial s_j }$ and $ Y = f^\varphi_{dd'} $ or $ g =  \frac{\partial^2 \psi_d}{\partial s_i \partial s_j} $ and $ Y = f^\varphi_{d} $ allows us to provide the desired bounds. Both possible choices of $ g $ are continuous on $ \phi\left( \cl\left( \mathcal{W}_s \right) \right)$ and hence $\max_{x \in \phi(\mathcal{W}_s)}\{\, g^2(x)\,\} < \infty $. A similar argument establishes finiteness for the maximum of the expected value of $Y$ squared. Furthermore both choices of $g$ are differentiable, as $\psi$ is $C^3$, and thus Lipschitz on $ \phi\left( \cl\left( \mathcal{W}_s \right) \right)$. Also,
	\begin{equation}
	\mathbb{E}\Bigg[\, \Bigg( f^\varphi_{dd'}\big(\psi(x)\big) - f^\varphi_{dd'}\big(\psi(y)\big) \Bigg)^2 \,\Bigg]
	\leq L \Vert\, \psi(x) - \psi(y) \,\Vert^{2\eta}
	\leq \tilde L \Vert\, x - y \,\Vert^{2\eta}
	\end{equation}
	for some $L, \tilde{L} \in (0, \infty)$ by \textit{\textbf{(C1)}} and Lipschitz continuity of $\psi$. A similar argument works for $ Y = f^\varphi_{d} $ since $f$ is $C^2$ and Gaussian.
\end{proof}

The key observation, that is needed to prove our main result, is the following extension of Lemma 4.2
from \cite{Cheng:2015extremes} to certain submanifolds of $ \tilde{\mathcal{M}} $.
\begin{lemma}\label{lem:SimpleSubmanifolds}
	Let $ \mathcal{N} \subset \tilde{\mathcal{M}} $ be a relatively compact and
	extendable $ d $-dimensional submanifold for some $ d \in \mathbb{N}. $ Let $ \tilde f $ be a Gaussian random field on $ \tilde{ \mathcal{M} } $ satisfying
	Assumption \ref{assumption:Fields} and let $ f $
	denote the restriction of $ \tilde f $ to $ \mathcal{N} $.
	Then $ \mE\big[\,  \mu_f\big( \mathcal{N} \big)^2 \,\big] < \infty $.
\end{lemma}
\begin{proof}
	Since $ \mathcal{N} $ is extendable there exists a $d$-dimensional
	submanifold $ \mathcal{A} $ without boundary such that
	$ \cl( \mathcal{N} ) \subset \mathcal{A} $.
	For $ s\in \cl( \mathcal{N} ) $, let $ (\,\mathcal{U}_s, \varphi \,) $
	be a chart around $ s $ satisfying Assumption \ref{assumption:Fields}.
	Let $ \left(\, \mathcal{U}^\mathcal{A}_s, \varphi_\mathcal{A} \,\right) $
	be a submanifold chart of $ \mathcal{A} $ around $ s $, i.e. such that,
	\begin{equation}\label{eq:submfChart}
	\varphi_\mathcal{A}\big( \mathcal{U}^\mathcal{A}_s  \cap
	\mathcal{A} \big) = \left\{\, (\,x_1,\ldots, x_D\,) \in \mathbb{R}^D: ~ x_{d+1} = \ldots = x_D = 0 \,\right\} \cap  \varphi_\mathcal{A}\big( \mathcal{U}^\mathcal{A}_s \big)\,.
	\end{equation}
	
	By Lemma 1, there exists some non-empty open set $\mathcal{W}_s \subset \mathcal{U}_s^{\mathcal{A}}$ such that $f^{\varphi_{\mathcal{A}}}$ satisfies assumptions \textbf{(C1)} and \textbf{(C2)} on $\mathcal{W}_s$. By \eqref{eq:submfChart} it follows that
	
	\begin{equation}
	\varphi_\mathcal{A}\big( \mathcal{W}_s  \cap
	\mathcal{A} \big) = \left\{\, (x_1,\ldots, x_D) \in \mathbb{R}^D: ~ x_{d+1} = \ldots = x_D = 0 \,\right\} \cap  \varphi_\mathcal{A}\big( \mathcal{W}_s \big)\,.
	\end{equation}
	
	Define the projection $\pi: \varphi_{\mathcal{A}}\big(\mathcal{W}_s \cap \mathcal{A}\big) \rightarrow \mathcal{V}_s, (\,x_1, \dots, x_D\,) \mapsto (\,x_1, \dots, x_d\,)$, where $ \mathcal{V}_s = \{\,  (\, e_1^T\varphi_\mathcal{A}(y),\ldots, e_d^T\varphi_\mathcal{A}(y)\,):~ y\in\mathcal{U}^\mathcal{A}_s  \cap \mathcal{A} \,\} \subset \mathbb{R}^d $ and $\lbrace e_i \rbrace_{i = 1, \dots, d}$ is the standard basis on $\mathbb{R}^d$. Then $\pi \circ \varphi_{\mathcal{A}}: \mathcal{W}_s \cap \mathcal{A}\rightarrow \mathcal{V}_s$ is a diffeomorphism since $\pi^{-1}:~\mathcal{V}_s \rightarrow \varphi_\mathcal{A}(\mathcal{W}_s)$ can be defined as the differentiable map $(x_1,\ldots,x_d)\mapsto(x_1,\ldots,x_d, 0,\ldots, 0)\in \mathbb{R}^D$. Now, $f^{\pi \circ \varphi_{\mathcal{A}}}=f^{\varphi_{\mathcal{A}}}\circ\pi^{-1}: \mathcal{V}_s \rightarrow \mathbb{R}$ satisfies conditions (C1) and (C2) from \cite{Cheng:2015extremes} since the same calculation as in equations \eqref{eq:1stDeriv} and \eqref{eq:2ndDeriv} shows that, for $x \in \mathcal{V}_s$, the random vectors obtained from $f^{\varphi_{\mathcal{A}}}\circ\pi^{-1}(x)$ in these conditions are identical to the random vectors obtained from $f^{\varphi_{\mathcal{A}}}$, evaluated at $\pi^{-1}(x)$, if the coordinates with indices $ i,j \in \{\, d+1,\ldots, D \,\} $ are removed. Therefore Lemma 4.2 from \cite{Cheng:2015extremes} can be applied and as such there exists $\epsilon(s) > 0$ such that $\mu_{f^{\pi \circ \varphi_\mathcal{A}}}\big(\,B_{\epsilon(s)}\big(\pi\circ\varphi_{\mathcal{A}}(s)\big)\,\big) < \infty$.
	
	For what follows, let $ \mathcal{O}_s =  \varphi_{\mathcal{A}}^{-1} \circ \pi^{-1} \big(\, B_{\epsilon(s)}\big(\pi \circ\varphi_{\mathcal{A}}(s)\big) \,\big)$. The family $\left\{\, \mathcal{O}_s \cap \mathcal{A}  \,\right\}_{s\in\cl( \mathcal{N}) } $ forms an open cover of the compact space $ \cl( \mathcal{N} ) $. Hence there exist $s_1,\ldots,s_K$, $K \in \mathbb{N}$, such that $ \bigcup_{k=1}^K \mathcal{O}_{s_k} \cap \mathcal{A} \supset \cl( \mathcal{N} ) $.
	Now,
	\begin{equation}
	\begin{split}
	\mu_{f}\big( \mathcal{N} \big) \leq \mu_{f}\big( \cl( \mathcal{N} ) \big) \leq \sum_{k=1}^K \mu_f\Big( \mathcal{O}_{s_k} \cap \mathcal{A} \Big)
	= \sum_{k=1}^K \mu_{f^{\pi \circ \varphi_{\mathcal{A}}}}\Big( B_{\epsilon(s_k)}\big(\pi \circ\varphi_{\mathcal{A}}(s)\big) \Big)        
	\end{split},
	\end{equation}
	since $\pi \circ \varphi_{\mathcal{A}}$ is a diffeomorphism, and so it follows that
	\begin{equation*}
	\begin{split}
	\mE\Big[\, \mu_f\big( \mathcal{N} \big)^2 \,\Big]
	\leq&
	\sum_{ k = 1 }^K \mE\Bigg[\, \mu_{f^{\pi \circ \varphi_{\mathcal{A}}}}\Big(B_{\epsilon(s_k)}\big(\pi \circ\varphi_{\mathcal{A}}(s)\big) \Big)^2 \,\Bigg] \\ ~ ~ ~ ~& +
	\sum_{ k = 1 }^K\sum_{ \substack{k' = 1\\ k'\neq k} }^K \mE\Bigg[\, \mu_{f^{\pi \circ \varphi_{\mathcal{A}}}}\Big( B_{\epsilon(s_k)}\big(\pi \circ\varphi_{\mathcal{A}}(s)\big) \Big)
	\mu_{f^{\pi \circ \varphi_{\mathcal{A}}}}\Big( B_{\epsilon(s_{k'})}\big(\pi \circ\varphi_{\mathcal{A}}(s)\big) \Big) \,\Bigg]\\
	\leq&
	\sum_{ k = 1 }^K \mE\Bigg[\, \mu_{f^{\pi \circ \varphi_{\mathcal{A}}}}\Big( B_{\epsilon(s_k)}\big(\pi \circ\varphi_{\mathcal{A}}(s)\big) \Big)^2 \,\Bigg] \\ ~ ~ ~ ~& +
	\sum_{ k = 1 }^K\sum_{ \substack{k' = 1\\ k'\neq k} }^K \sqrt{ \mE\Bigg[\, \mu_{f^{\pi \circ \varphi_{\mathcal{A}}}}\Big( B_{\epsilon(s_k)}\big(\pi \circ\varphi_{\mathcal{A}}(s)\big) \Big)^2 \,\Bigg]}\sqrt{ \mE\Bigg[\, \mu_{f^{\pi \circ \varphi_{\mathcal{A}}}}\Big( B_{\epsilon(s_{k'}t)}\big(\pi \circ\varphi_{\mathcal{A}}(s)\big) \Big)^2 \,\Bigg]}\,.
	\end{split}
	\end{equation*}
	The right hand side of this inequality is finite since all
	expectations in the sums are finite by Lemma 4.2 from \cite{Cheng:2015extremes}. 
\end{proof}

\begin{theorem}
	Let $ \mathcal{M} = \bigsqcup_{ d = 0 }^D \pdM \subset \tilde{\mathcal{M}} $,
	where $ \pdM $, for $ d = \{\, 1, \ldots,D \,\} $, is either empty or a disjoint union of finitely
	many $ d $-dimensional relatively compact and extendable $ C^3 $-submanifolds of
	$ \tilde{\mathcal{M}} $ without boundary and $ {\partial_0 \cM} $ is either empty or
	is the union of finitely many points of $ \tilde{\mathcal{M}} $.
	Let $ \tilde f $ be a Gaussian random field on $ \tilde{ \mathcal{M} } $ satisfying
	Assumption \ref{assumption:Fields} and let $ f $ denote its restriction to $ \mathcal{M} $.
	Then $ \mE\Big[  \mu_f^2\big( \mathcal{M} \big) \Big] < \infty $.
\end{theorem}
\begin{proof}
	The number of critical points of $ \mathcal{M} $ is given by
	$ \mu_f\big(\mathcal{M}\big) = \sum_{d=1}^D \mu_{f\vert_{\partial_d\mathcal{M}}}\big( \partial_d\mathcal{M}\big) $. Thus, applying Lemma \ref{lem:SimpleSubmanifolds} to $\partial_d \mathcal{M}$ and the Cauchy-Schwartz inequality it follows that 
	\begin{equation*}
	\begin{split}
	\mE\Big[\, \mu_f^2\big(\mathcal{M}\big) \,\Big]
	= \sum_{d=1}^D\sum_{d'=1}^D  \mE\Big[\, \mu_{f\vert_{\partial_d\mathcal{M}}}\big( \partial_d\mathcal{M}\big) \mu_{f\vert_{\partial_d\mathcal{M}}}\big( \partial_{d'}\mathcal{M}\big) \,\Big]
	\leq \left(\, \sum_{d=1}^D  \sqrt{ \mE\Big[\, \mu_{f\vert_{\partial_d\mathcal{M}}}^2\big( \partial_d\mathcal{M}\big) \,\Big] }\,\right)^2
	< \infty\,.
	\end{split}
	\end{equation*}
\end{proof}

%% file: appendix.tex
\appendix
\appendix
\section{Technical Lemmas}
\begin{lemma}\label{lem:Geman}
	Suppose that $\tilde{\mathcal{M}} \subset \mathbb{R}^D$, $D \in \mathbb{N}$, has non-empty interior and that $f: \tilde{\mathcal{M}} \rightarrow \mathbb{R}$ is a stationary Gaussian random field. Then if $f$ satisfies \textbf{(C1)} (with respect to the identity chart) it satisfies the generalized Geman condition on $\tilde{\mathcal{M}}$.
\end{lemma}
\begin{proof}
	Without loss of generality, assume that $ 0 $ lies within the interior of $\tilde{\mathcal{M}} $ and choose $ \delta > 0 $ such that $ B_\delta(0) \subset \tilde{\mathcal{M}} $. Let $r$ denote the covariance function of $f$. Then for all $ t \in B_\delta(0) $ and each $ 1 \leq i,j,k,l \leq D, $ letting $r_{ijkl} = \frac{\partial^4 r}{\partial t_i\partial t_j\partial t_k\partial t_l}$, 
	\begin{align*}
	\left\vert\, r_{ijkl}(0) - r_{ijkl}(t) \,\right\vert &= 
	\left|\, \mathbb{E}\left[\, f_{ij}(0)f_{kl}(0) \,\right] - \mathbb{E}\left[\, f_{ij}(0)f_{kl}(t) \,\right]\,\right| 
	= \left|\, \mathbb{E}\big[\, f_{ij}(0)\,\big(f_{kl}(0) -f_{kl}(t)\big) \,\big]\,\right| \\
	&\leq \sqrt{\,\mathbb{E}\left[\, f_{ij}(0)^2 \,\right]\,}\sqrt{\, \mathbb{E}\left[\, \big(\,f_{kl}(0) - f_{kl}(t)\,\big)^2 \,\right]\,} \leq 
	M L^{1/2}\left\lVert t \right\rVert^{\eta}.
	\end{align*}
	where $ M = \max_{1 \leq i,j \leq D}\sqrt{\,\mathbb{E}\left[ f_{ij}(0)^2 \right]\,} $ is finite because $ f_{ij} $ is a Gaussian random field. In particular, since $ \eta - D > -D, $
	\begin{equation*}
	\int_{\left\lVert t \right\rVert < \delta} \frac{\left\lVert\, r^{(4)}(0) - r^{(4)}(t) \,\right\rVert}{\left\lVert t \right\rVert^{D}}
	\leq
	\sum_{1 \leq i,j \leq D} \int_{\left\lVert t \right\rVert < \delta} \frac{ \left|\, r_{ijkl}(0) - r_{ijkl}(t) \,\right|}{\left\lVert t \right\rVert^{D}} \leq
	D^2	M L^{1/2} \int_{\left\lVert t \right\rVert < \delta} \left\lVert t \right\rVert^{\eta - D}
	<
	\infty.
	\end{equation*}
\end{proof}

\begin{lemma}\label{lem:MU}
	Let  $ U  \in \mathbb{R}^{D} $ be a random vector. Then for all $M_1, M_2 \in \mathbb{R}^{D \times D}$ such that $\det(M_1)\det(M_2)=1$ we have that
	\begin{equation}
		{\rm detcov }\big[ U \big] = \det\!\Big( { \rm Cov }\big[\, U, U \,\big] \Big)
			= \det\!\Big( { \rm Cov }\big[ \,M_1U, M_2U \,\big] \Big)
	\end{equation}
	In particular, for any invertible $M \in \mathbb{R}^{D \times D}$, this implies that
	\begin{equation}
		{\rm detcov }\big[ U \big]
			= {\rm detcov }\big[ M U \big] \,/ 
           \det\big( M \big)^{2D}\,.
	\end{equation}
\end{lemma}
\begin{proof}
    The determinant multiplication theorem  yields
    \begin{equation*}
        \begin{split}
            {\rm detcov }\big[ U \big]
            = \det(M_1){\rm detcov }\big[ U \big] \det(M_2^T)
           = \det\!\Big( M_1{ \rm Cov }\big[ U \big] M_2^T\Big)
           = \det\!\Big( { \rm Cov }\big[\, M_1U, M_2U \,\big] \Big)\,.
        \end{split}
    \end{equation*}
    Taking $ M_1 = M_2 = M\,/\det(M)$ establishes the second claim.
\end{proof}